\documentclass[12pt]{article}
\usepackage[cp1251]{inputenc}
\usepackage[english]{babel}
\usepackage{amssymb,amsmath}
 \usepackage{color}
\newtheorem{Def}{Definition}[section]
\newtheorem{Ex}{Example}[section]
\newtheorem{Th}{Theorem}[section]

\newtheorem{Prop}{Proposition}[section]
\newtheorem{Cor}{Corollary}[section]

\newenvironment{Proof}
{\par\noindent{\bf Proof.}}
{\hfill$\scriptstyle\blacksquare$}

\title{On a metric on the space of monetary risk measures}
\author{Sh.~A.~Ayupov, A.~A.~Zaitov}

\begin{document}

\maketitle
\thispagestyle{empty}

\begin{center}

\end{center}
\begin{abstract}
We introduce a metric on the space of monetary risk measure, which generates the point-wise convergence topology and extends the metric on the initial compactum.\\


\end{abstract}

\tableofcontents

\section{Introduction}

The financial market faces risks arising from many types of uncertain losses, including market risk, credit risk, liquidity risk, operational risk, etc. In 1988, the Basel Committee on Banking Supervision proposed measures to control credit risk
in banking. A risk measure called the \textit{value-at-risk}, acronym VaR, became, in the 1990s, an important tool of risk assessment and management for banks, securities companies, investment funds, and other financial institutions in asset allocation and
performance evaluation. The VaR associated with a given confidence level for a venture capital is the upper limit of possible losses in the next certain period of time. In 1996 the Basel Committee on Banking Supervision endorsed the VaR as one of the acceptable methods for the bank’s internal risk measure. However, due to the defects of VaR, a variety of new risk measures came into being. This paper
focuses on the definition the distance between normed monetary risk measures. For an overview of
the subject we refer to \cite{Yan2018}.

Assume that all possible states and events that may occur at the terminal time are known, namely, a measurable space $(X,\, \mathcal{F})$ is given. The financial position (here refers to the wealth deducted investment cost) is usually described by a measurable
function $\varphi$ on $(X,\, \mathcal{F})$. If we assume that a probability measure $\mathbb{P}$ is given on measurable space $(X,\, \mathcal{F})$, the financial positions is usually described by a random variables $\widetilde{\varphi}$ on $(X,\, \mathcal{F},\, \mathbb{P})$. In order to facilitate the notation and description, we use $\varphi = - \widetilde{\varphi}$ to denote the \textit{potential loss} at the terminal time of trading. Here the potential loss is relative to a reference point in terms. If $\varphi$ is a negative value, it indicates a \textit{surplus}.

A risk measure is a numerical value $\mu(\varphi)$ to quantify the risk of a financial position (it may be a potential loss or a surplus as well) $\varphi$. If we denote the set of financial positions to be considered by $\mathcal{G}$, a risk measure $\mu$ is a map from $\mathcal{G}$ to $\mathbb{R}$. Usually, they take $L^{\infty}(X,\, \mathcal{F},\, \mathbb{P})$ or $L^{\infty}(X,\, \mathcal{F})$ as the set of all financial positions $\mathcal{G}$ or $\mathcal{G}(\mathbb{R})$, where the former is the set of all bounded $\mathcal{F}$ measurable functions on $(X,\, \mathcal{F})$, endowed with the uniform norm $\|\cdot\|_{\infty}$, and the latter is the set of equivalence classes of the former under probability $\mathbb{P}$. In the former case, the states and the probabilities of the possible events are unknown or are not consensus in the market, and then the risk measure is called model-free. In the latter case, the risk measure is called model-dependent. In the model-dependent
case, naturally, we always assume that the risk measure $\mu$ satisfies the following property: If $\varphi = \psi$, $\mathbb{P}$-a.s., then $\mu(\varphi) = \mu(\psi)$.

\begin{Def}\label{monetaryriskL}
{\rm A map $\mu$ from $\mathcal{G}$ to $\mathbb{R}$ is called a \textit{monetary risk measure}, abbreviated
as \textit{risk measure}, if it satisfies two conditions:
\begin{itemize}
\item[$(1)$] \textit{monotonicity}: For all $\varphi,\, \psi \in \mathcal{G}$ satisfying $\varphi\le \psi$, it holds that $\mu(\varphi)\le \mu(\psi)$.
\item[$(2)$] \textit{translation invariance}: For all $\varphi \in \mathcal{G}$ and any real number $\alpha$, it holds that
\begin{gather*}
\mu(\varphi + \alpha) = \mu(\varphi) + \alpha.
\end{gather*}
\end{itemize}}
\end{Def}

It is known that the algebra $L^{\infty}(X,\, \mathcal{B}(X),\, \mathbb{P})$ is isomorphic to the algebra  $C(X)$ of all continuous functions on a compact  $X$ (to be more precise, $X$ is so called hyperstonean compact). Thus the above notion of financial position  can be interpreted as an element of the algebra $C(X)$,  while the monetary risk measure can be considered as a map from $C(X)$ to   $\mathbb{R}$.

In the present paper we consider the problem in a more general setting, where financial positions are interpreted as elements of the algebra $C(X)$, where $X$ is a compact Hausdorff space. Then the algebra $C(X)$ plays a role of the set $\mathcal{G}$ of all financial positions. A \textit{risk measure} is a numerical value $\mu(\varphi)$ to quantify the risk of a financial position $\varphi\in C(X)$. A map $\mu$ from $C(X)$ to $\mathbb{R}$ is called a \textit{monetary risk measure} on $X$, if it satisfies the conditions $(1)$ and $(2)$ of Definition \ref{monetaryriskL}.

A monetary risk measure $\mu\colon C(X)\to \mathbb{R}$ is called \textit{normed} if
\begin{itemize}
\item[$(3)$] $\mu(1_X)=1$.
\end{itemize}

In works  \cite{Rad1998} -- \cite{Zait2019stemm} normed monetary risk measures are called as order-preserving functionals and the set of such functionals is denoted by $O(X)$. The mentioned papers were devoted to study of $O(X)$. We consider $O(X)$ as a subspace of the Tychonoff product $\mathbb{R}^{C(X)}$. The base of the induced topology consists of the sets of the form
\begin{gather*}
\langle\mu;\, \varphi_{1},\, \dots,\, \varphi_{n};\, \varepsilon\rangle = \{\nu\in O(X):\, |\nu(\varphi_{i}) - \mu(\varphi_{i})| < \varepsilon,\, i=1,\, \dots,\, n\},
\end{gather*}
where $\mu\in O(X)$, $\varphi_{i}\in C(X)$, $i=1,\, \dots,\, n$, and $\varepsilon > 0$. Note that the induced topology and point-wise convergence topology coincide. For every compact Hausdorff space $X$ the space $O(X)$ is also a compact Hausdorff space. $O(X)$ is a compact sublattice of $\mathbb{R}^{C(X)}$.

A closed subset $A$ of $O(X)$ is called $O$-\textit{convex}  if for each $\mu\in O(X)$ with $\inf\, A \le \mu \le \sup\, A$ we have $\mu\in A$. It is known (\cite{Rad1999}, Lemma 3) that for each map $f\colon X\to Y$ and $\nu\in O(Y)$ the preimage $O(f)^{-1}(\nu)$ is an $O$-convex subset of $O(X)$. \label{preimageOconvex}

\begin{Prop}\label{interOconvex}
{\rm Let $A_{1}$, $A_{2}$ be $O$-convex subsets of $O(X)$. Then $A_{1}\cap A_{2}$ is an $O$-convex subset of $O(X)$.}
\end{Prop}

\begin{Proof} The proof consists of directly checking.

\end{Proof}

Let $X$, $Y$ be compact Hausdorff spaces, $f\colon X\to Y$ be a continuous map. Then a map $O(f)\colon O(X)\to O(Y)$, defined as $O(f)(\mu)(\varphi) = \mu(\varphi\circ f)$, $\varphi\in C(Y)$, is continuous. So, the monetary risk measures $\mu$ and $O(f)(\mu)$ act \textbf{the same rule}. Let us denote $\mu_{X} = \mu$, $\mu_{Y} = O(f)(\mu)$. Then the discussed situation is bring to light: $\mu_{X}(\psi) = \mu(\psi)$, $\psi\in C(X)$, and $\mu_{Y}(\varphi) = \mu(\varphi\circ f)$, $\varphi\in C(Y)$, i.~e. $\mu_{X}$ and $\mu_{Y}$ both act by the rule $\mu$.

Let $\mu_{i}\in O(X)$, $i=1,\, 2$. We say that $\mu_{1} = \mu_{2}$ if $\mu_{1}(\varphi) = \mu_{2}(\varphi)$ for all $\varphi\in C(X)$. The following statement is rather obvious.

\begin{Prop}\label{meas.equiv}
{\rm Let $\mu_{i}\in O(X)$, $i=1,\, 2$. Then $\mu_{1} = \mu_{2}$ if and only if $\mbox{supp}\, \mu_{1} = \mbox{supp}\, \mu_{2}$ and $\mu_{2} = O(h)(\mu_{1})$, where $h\colon\, \mbox{supp}\, \mu_{1} \rightarrow \mbox{supp}\, \mu_{2}$ is the identity map}.
\end{Prop}

Note that if $f$ is a surjective map then $O(f)$ is also a surjective map. If $X$ is a closed subset of $Y$ and $f$ is an embedding then $O(f)$ is also an embedding.

Let $X$ be a compact Hausdorff space and $\mu$ be a monetary risk measure. $\mu$ is \textit{concentrated} on a closed subset $A$ of $X$ if  $\mu\in O(A)$. Note that (Lemma 4, \cite{Rad1998}) for a closed subset $A\subset X$ a monetary risk measure $\mu\in O(X)$ is supported on $A$ if and only if for every pair $\varphi,\, \psi\in C(X)$ such that $\varphi|_{A} = \psi|_{A}$ one has $\mu(\varphi) = \mu(\psi)$. \label{coinsidemeasures} The smallest (with respect to inclusion) closed subset $\mbox{supp}\mu$ of $X$ on which $\mu$ is concentrated is said to be a \textit{support} of monetary risk measure $\mu$. Evidently,
\begin{gather*}
\mbox{supp}\,\mu = \cap\left\{A:\, A\, \mbox{ is a closed set in }\, X\, \mbox{ and }\, \mu\in O(A)\right\}.
\end{gather*}

For a point $x\in X$ the Dirac measure $\delta_{x}$, defined by $\delta_{x}(\varphi) = \varphi(x)$, $\varphi\in C(X)$, is a monetary risk measure, concentrated at the singleton $\{x\}$, i.~e. $\mbox{supp}\,\delta_{x} = \{x\}$.

A subset $L\subset C(X)$ is called an $A$-\textit{subspace} if $0_{X}\in L$ and for every $\varphi\in L$ and every $c\in \mathbb{R}$ we have $\varphi+c_{X}\in L$. According to the analog of the Hahn-Banach theorem (see, \cite{Rad1998}, \cite{Zait2005mfat}) for every normed monetary risk measure $\mu \colon\, L\to \mathbb{R}$ there exists a normed monetary risk measure $\widetilde{\mu}\colon\, C(X)\to \mathbb{R}$ such that $\widetilde{\mu}|_{L} = \mu$. \label{HahnBanach}\\

The space $O(X)$ of monetary risk measure does not embed into any linear space with finite algebraic dimension if $X$ consists more than one point.

\begin{Ex}{\rm \cite{Zait2019stemm}.
Let $X=\{0,\, 1\}$ be a discrete two-point space. Then $C(X) = \mathbb{R}^2$. Each functional $\mu\colon C(X)\to \mathbb{R}$ defined  by the equality
\begin{multline}\label{exam}
\mu(\varphi) = \alpha_1\,\varphi(0) + \alpha_2\,\varphi(1) + \alpha_3\,\max\{\varphi(0)+\lambda_1,\, \varphi(1)+\lambda_2\} + \\
+ \alpha_4\,\min\{\varphi(0)+\lambda_3,\, \varphi(1)+\lambda_4\} + \alpha(\varphi)\,f(\varphi(1)-\varphi(0))
\end{multline}
is a normed monetary risk measure. Here $\alpha_1 + \alpha_2 + \alpha_3 + \alpha_4=1$ with $\alpha_i\ge 0$, $i=1,\, 2,\, 3,\, 4$, $\lambda_1,\, \lambda_2\in [-\infty,\, 0]$ with $\max\{\lambda_1,\, \lambda_2\}=0$,  $\lambda_3,\, \lambda_4\in [0,\, +\infty]$ with $\min\{\lambda_3,\, \lambda_4\}=0$,
\begin{equation}
\alpha(\varphi)
=\begin{cases}
\min\{\alpha_1,\, \alpha_2\}, & \mbox{ if }  \alpha_3 = \alpha_4 =0,\\
\min\{\alpha_1+\alpha_3+\alpha_4,\, \alpha_2\}, & \mbox{ if }  \varphi(0)+\lambda_1 \ge \varphi(1)+\lambda_2 \mbox{ and } \varphi(0)+\lambda_3 \le \varphi(1)+\lambda_4,\\
\min\{\alpha_1+\alpha_3,\, \alpha_2+\alpha_4\}, & \mbox{ if }  \varphi(0)+\lambda_1 \ge \varphi(1)+\lambda_2 \mbox{ and } \varphi(0)+\lambda_3 > \varphi(1)+\lambda_4,\\
\min\{\alpha_1+\alpha_4,\, \alpha_2+\alpha_3\}, & \mbox{ if }  \varphi(0)+\lambda_1 < \varphi(1)+\lambda_2 \mbox{ and } \varphi(0)+\lambda_3 \ge \varphi(1)+\lambda_4,\\
\min\{\alpha_1,\, \alpha_2+\alpha_3+\alpha_4\}, & \mbox{ if }  \varphi(0)+\lambda_1 < \varphi(1)+\lambda_2 \mbox{ and } \varphi(0)+\lambda_3 > \varphi(1)+\lambda_4,
\end{cases}
\label{examfunc}
\end{equation}
and, finally, $f\colon \mathbb{R}\to \mathbb{R}$ is a continuous non-decreasing function such that
\begin{itemize}
\item[$(1^*)$]  $f(0) = 0$,
\item[$(2^*)$] $t \le f(t) \le 0$ and concave at $t \le 0$,
\item[$(3^*)$] $0 \le f(t) \le t$ and convex at  $t \ge 0$.
\end{itemize}

Since the set of functions $f$ considered in (\ref{exam}) and satisfying conditions $(1^*) - (3^*)$, is uncountable, it follows that the space $O(X)$ of normed monetary risk measure can not be embedded in any space with finite (even countable) algebraic dimension if the compact $X$ contains more than one point.}
\end{Ex}

Example \ref{exam} shows that there exists extremely many monetary risk measures in practice. Further, a question arises wether one can evaluate the difference between distinct monetary risk measures. In other words, \textbf{is it possible to specify distance between monetary risk measures?}

It is known \cite{Rad1998} that for a compact Hausdorff space $X$ the equality $w(X) = w(O(X))$ holds, where $w$ is the weight (i.~e. the smallest cardinal number which is the cardinality of an open base) of a topological space. From here follows that if $X$ is a compactum then $O(X)$ is a compactum, i.~e. the space of normed monetary risk measures is metrisable compact space. Though for a metrisable compact space $X$ the space $O(X)$ is metrisable, but the rule of definition distance between monetary measures still was not known. In the present paper for a given compact metric space $(X,\, \rho)$ we introduce a metric $\rho_{O}$ on the space $O(X)$ of normed monetary risk measure, which generates the point-wise convergence topology on $O(X)$. Besides, $\rho_{O}$ is an extension on $O(X)$ of the metric $\rho$.

\section{Basic constructions}

For a compact Hausdorff space $X$ we put
\begin{gather*}
X_1 = X_2 = X_3 = X,\qquad X_{1\,2\,3} = X^3 = X_1\times X_2\times X_3,\qquad X_{i\,j} = X^2 = X_i\times X_j,
\end{gather*}
and let
\begin{gather*}
\pi^{1\,2\,3}_{i\,j}\colon X_{1\,2\,3}\to X_{i\,j},\qquad  \pi^{i\,j}_{k}\colon X_{i\,j}\to X_k,\qquad 1\le i< j\le 3,\qquad k\in \{i,\,j\},
\end{gather*}
be the corresponding projections.

Obviously, that
\begin{gather*}
O(\pi^{i\,j}_{i})^{-1}O(X_{i})\cap O(\pi^{i\,j}_{j})^{-1}O(X_{j}) = O(X_{i\,j}),\qquad 1\le i< j\le 3,\\
\bigcap\limits_{1\le i< j\le 3}O(\pi^{1\,2\,3}_{i\,j})^{-1}O(X_{i\,j}) = O(X_{1\,2\,3}).
\end{gather*}

\begin{Th}\label{nonemptyofpreimage}
{\rm For every pair $(\mu_{1},\, \mu_{2})\in O(X_1)\times O(X_2)$ we have
\begin{gather*}
O(\pi^{1\,2}_{1})^{-1}(\mu_{1})\cap O(\pi^{1\,2}_{2})^{-1}(\mu_{2})\ne \varnothing.
\end{gather*}}
\end{Th}

\begin{Proof}
At first we consider a particular case: let $\mu_{1} = \delta_{x}$ and $\mu_{2} = \delta_{y}$, where $x\in X_{1}$, $y\in X_{2}$. The spaces $O(X_{1}\times \{y\})$ and $O(X_{1})\times \{\delta_{y}\}$ are homeomorphic. Indeed, one may determine the homeomorphism as the correspondence $O(X_{1}\times \{y\})\ni \lambda_{y}\mapsto (\lambda,\, \delta_{y})\in O(X_{1})\times \{\delta_{y}\}$, where $\mbox{supp}\,\lambda_{y} = \{(x,\, y):\, x\in \mbox{supp}\,\lambda\}$, and monetary risk measures $\lambda_{y}$ and $\lambda$ act the same rule, i.~e.  $\lambda_{y}(\varphi) = \lambda(\varphi\circ \pi^{1\,2}_{1})$, $\varphi\in C(X_{1}\times X_{2})$. Consequently, $O(\pi^{1\,2}_{2})^{-1}(\delta_{y}) = O(X_{1}\times\{y\})$. Similarly, $O(\pi^{1\,2}_{1})^{-1}(\delta_{x}) = O(\{x\}\times X_{2})$. It is easy to see that $\delta_{(x,\,y)}$ is a unique monetary risk measure such that $\delta_{(x,\,y)}\in O(\{x\}\times X_{2})\cap O(X_{1}\times\{y\})$. Thus, $O(\pi^{1\,2}_{1})^{-1}(\delta_{x})\cap O(\pi^{1\,2}_{2})^{-1}(\delta_{y})\ne \varnothing$.

Now we consider the general case. We construct a set
\begin{gather*}
\mathfrak{B} = \{\varphi\circ \pi^{1\,2}_{i} + c_{X_{1}\times X_{2}}:\, \varphi\in C(X),\, c\in \mathbb{R},\, i = 1,\, 2\}.
\end{gather*}
Then $\mathfrak{B}$ is an $A$-subspace in $C(X_{1}\times X_{2})$. Define a functional $\widetilde{\mu}_{1\, 2}\colon\mathfrak{B}\to \mathbb{R}$ as following
\begin{gather*}
\widetilde{\mu}_{1\, 2}(\varphi\circ \pi^{1\,2}_{i} + c_{X_{1}\times X_{2}}) = \mu_{i} (\varphi) + c,\qquad i = 1,\, 2.
\end{gather*}
It is easy to see that $\widetilde{\mu}_{1\, 2}$ is translation invariance and normed.

We will show $\widetilde{\mu}_{1\, 2}$ has monotonicity property. Taking $\varphi\circ \pi^{1\,2}_{i}$, $\psi\circ \pi^{1\,2}_{i}$ $\in \mathfrak{B}$ with $\varphi\circ \pi^{1\,2}_{i} \le \psi\circ \pi^{1\,2}_{i}$, we obtain $\varphi \le \psi$, and thence
\begin{gather*}
\widetilde{\mu}_{1\, 2} (\varphi\circ \pi^{1\,2}_{i}) = \mu_{i}(\varphi) \le \mu_{i}(\psi) = \widetilde{\mu}_{1\, 2} (\psi\circ \pi^{1\,2}_{i}),\qquad i = 1,\, 2.
\end{gather*}
Take now $\varphi\circ \pi^{1\,2}_{1}$, $\psi\circ \pi^{1\,2}_{2}$ $\in \mathfrak{B}$ such, say, that $\varphi\circ \pi^{1\,2}_{1} \le \psi\circ \pi^{1\,2}_{2}$. Then
\begin{gather*}
\max\{\varphi(x):\, x\in X\}\le \min\{\psi(x):\, x\in X\}.
\end{gather*}
Choosing any $a\in [\max\{\varphi(x):\, x\in X\},\, \min\{\psi(x):\, x\in X\}]$, we see
\begin{gather*}
\widetilde{\mu}_{1\, 2} (\varphi\circ \pi^{1\,2}_{1}) = \mu_{1}(\varphi) \le a \le \mu_{2}(\psi) = \widetilde{\mu}_{1\, 2} (\psi\circ \pi^{1\,2}_{2}).
\end{gather*}
One similarly can establish the monotonicity of $\widetilde{\mu}_{1\, 2}$ in the case when $\varphi\circ \pi^{1\,2}_{1} \ge \psi\circ \pi^{1\,2}_{2}$.

Thus, $\widetilde{\mu}_{1\, 2}$ is a normed monetary risk measure on the $A$-subpace $\mathfrak{B}$. By the analog of the Hahn-Banach theorem (see Page \pageref{HahnBanach} of the present paper)  $\widetilde{\mu}_{1\, 2}$ has an extension $\mu_{1\, 2}$ all over $C(X_{1}\times X_{2})$, which is a normed monetary risk measure. We have
\begin{gather*}
O(\pi^{1\, 2}_{i})(\mu_{1\, 2}) = \mu_{1\, 2}(\varphi\circ \pi^{1\, 2}_{i}) = \widetilde{\mu}_{1\, 2}(\varphi\circ \pi^{1\, 2}_{i}) = \mu_{i}(\varphi),\qquad i = 1,\, 2.
\end{gather*}
Consequently, $\mu_{1\, 2}\in O(\pi^{1\,2}_{1})^{-1}(\mu_{1})\cap O(\pi^{1\,2}_{2})^{-1}(\mu_{2})$.

\end{Proof}

Denote $\Lambda(\mu_{1},\, \mu_{2}) = O(\pi^{1\,2}_{1})^{-1}(\mu_{1})\cap O(\pi^{1\,2}_{2})^{-1}(\mu_{2})$. An element $\xi\in \Lambda(\mu_{1},\, \mu_{2})$ we call as a $(\mu_{1},\, \mu_{2})$-\textit{admissible} monetary risk measure.

\begin{Cor}\label{overmu2}
{\rm
For normed monetary risk measures
\begin{gather*}
\mu_2 \in O(X_2), \qquad \mu_{1\,2}\in O(X_{1\,2}),\qquad \mu_{2\,3}\in O(X_{2\,3})
\end{gather*}
such that
\begin{gather*}
O(\pi^{1\,2}_{2})(\mu_{1\,2}) = \mu_2 = O(\pi^{2\,3}_{2})(\mu_{2\,3}),
\end{gather*}
there exists a $\mu_{1\,2\,3}\in O(X_{1\,2\,3})$ satisfying the equalities
\begin{gather*}
O(\pi^{1\,2\,3}_{1\,2})(\mu_{1\,2\,3}) = \mu_{1\,2}\qquad \mbox{ and }\qquad O(\pi^{1\,2\,3}_{2\,3})(\mu_{1\,2\,3}) = \mu_{2\,3}.
\end{gather*}}
\end{Cor}

Really, to adopt this statement it is sufficient to note that
\begin{gather*}
O(\pi^{1\,2\,3}_{1\,2})^{-1}(\mu_{1\,2}) \cap O(\pi^{1\,2\,3}_{2\,3})^{-1}(\mu_{2\,3}) \ne \varnothing.
\end{gather*}
In this case we construct an $A$-subspace
\begin{gather*}
\mathfrak{B} = \{\varphi\circ \pi^{1\,2\,3}_{i\,(i+1)} + c_{X^{3}}:\, \varphi\in C(X^2),\, c\in \mathbb{R},\, i=1,\,2\}
\end{gather*}
in $C(X^{3})$ and repeat the analogous procedure as in the proof of Theorem \ref{nonemptyofpreimage}.

\begin{Prop}\label{diagneigh}
{\rm Let $X$ be a compactum and a sequence $\{\mu_{n}\}\subset O(X)$ converges to $\mu_{0}\in O(X)$ with respect to point-wise convergence topology. Then for every open neighbourhood $U$ of the diagonal $\Delta(X)=\{(x,\, x):\, x\in X\}$ of $X^{2}$ there exists a positive integer $n$ and for each $n'\ge n$ there exists a $(\mu_{0},\, \mu_{n'})$-admissible monetary risk measure $\mu_{0\, n'}\in O(X^2)$ such that $\mbox{supp}\, \mu_{0\, n'}\subset U$.}
\end{Prop}

\begin{Proof} The condition gives a sequence $\{\mbox{supp}\, \mu_{n}\}$ of closed subsets of $X$. It is well known that for a compact Hausdorff space $X$ its hyperspace $\exp\, X$ also is a compact Hausdorff space as well. Where $\exp\, X$ is the set of all nonempty closed subsets of $X$, and $\exp\, X$ is equipped with the Vietoris topology. So, the sequence $\{\mbox{supp}\, \mu_{n}\}$ has a limit. Let $A = \lim\limits_{n\to\infty}\mbox{supp}\, \mu_{n}$. Suppose $A \ne \mbox{supp}\, \mu_{0}$. Then Proposition \ref{meas.equiv} implies that $\mu_{0} \ne \lim\limits_{n\to\infty}\mu_{n}$. The got contradiction shows that $\lim\limits_{n\to\infty}\mbox{supp}\, \mu_{n} = \mbox{supp}\, \mu_{0}$.

Consider open neighbourhoods $V_{x}$ of points $x\in \mbox{supp}\, \mu_{0}$ such that $V_{x}\times V_{x} \subset U$. Since $\mbox{supp}\, \mu_{0}$ is a compact set, its open cover $\{V_{x}:\, x\in \mbox{supp}\, \mu_{0}\}$ has a finite subcover $\{V_{k}:\, k=1,\, \dots,\, l\}$, where $V_{k} = V_{x_{k}}$. Owing to convergence $\{\mbox{supp}\, \mu_{n}\}$ to $\mbox{supp}\, \mu_{0}$ there exists a positive integer $n$ such that $\mbox{supp}\, \mu_{n'} \in \langle V_{1},\, \dots,\, V_{l}\rangle$ for every $n' \ge n$. Here
\begin{gather*}
\langle V_{1},\, \dots,\, V_{l}\rangle = \left\{F\in \exp\,X:\, F\subset \bigcup\limits_{k=1}^{l}V_{k}\, \mbox{ and }\, F\cap V_{k}\ne \varnothing\, \mbox{ for each }\, k=1,\, \dots,\, l\right\}
\end{gather*}
is a basic open neighbourhood of $\mbox{supp}\, \mu_{0}$ with respect to the Vietoris topology in $\exp\, X$.

It is easy to see that $\mbox{supp}\, \mu_{n'} \in \langle V_{1},\, \dots,\, V_{l}\rangle$ if and only if
\begin{multline}\label{Hausd.top}
\mbox{supp}\, \mu_{n'} \subset \bigcup\limits_{k=1}^{l}V_{k} \quad \mbox{and for every}\quad  x\in \mbox{supp}\, \mu_{0}\quad
\mbox{there exists}\quad y\in\mbox{supp}\, \mu_{n'}\\
\quad \mbox{such that}\quad (x,\, y)\in V_{k}\times V_{k}\quad \mbox{for some}\quad k\in \{1,\, \dots,\, l\}.
\end{multline}

It remains to show that for every $n' \ge n$ there exists $\mu_{0\, n'}\in \Lambda(\mu_{0},\, \mu_{n'})$ such that  $\mbox{supp}\, \mu_{0\, n'} \subset \bigcup\limits_{k=1}^{l} V_{k}\times V_{k}$. Suppose, that for some $n'\ge n$ and for every $\xi\in \Lambda(\mu_{0},\, \mu_{n'})$ the intersection $\mbox{supp}\, \xi \cap (X^{2}\setminus U)$ is nonempty. Let $(x,\, y) \in \mbox{supp}\, \xi \cap (X^{2}\setminus U)$. Then $(x,\, y)\not\in V_{k}\times V_{k}$ for all $k=1,\, \dots,\, l$, which contradicts (\ref{Hausd.top}). The received contradiction finished the proof.

\end{Proof}

\section{On a metric on the space of monetary risk measures}

Let $(X,\, \rho)$ be a metric compact space. We suggest a distance function $\rho_{O}\colon O(X)\times O(X) \to \mathbb{R}$ as follows
\begin{equation}\label{rhoO}
\rho_{O}(\mu_1,\, \mu_2) = \inf\, \{\max\{\rho(x,\, y):\, (x,\, y)\in \mbox{supp}\, \xi\}:\, \xi\in \Lambda(\mu_{1},\, \mu_{2})\}.
\end{equation}

\begin{Prop}
{\rm For every pair $\mu_1$, $\mu_2$ of monetary risk measures there exists a $(\mu_{1},\, \mu_{2})$-admissible monetary risk measure $\mu_{1\, 2}\in O(X^{2})$ such that
\begin{gather*}
\rho_{O}(\mu_1,\, \mu_2) = \max\{\rho(x,\, y):\, (x,\, y)\in \mbox{supp}\, \mu_{1\, 2}\}.
\end{gather*}}
\end{Prop}

\begin{Proof} The Proof leans on Proposition \ref{interOconvex} and Lemma 3 \cite{Rad1999} (see Page \pageref{preimageOconvex} of the present paper).

\end{Proof}

\begin{Th}\label{met}{\rm
The function $\rho_{O}$ is a metric on the space $O(X)$ of monetary risk measures, which is an extension of the metric $\rho$ on $X$.}
\end{Th}

\begin{Proof} Since each $\xi\in O(X^2)$ is monotone then the inequality  $\rho\ge 0$ immediately implies $\rho_{O}\ge 0$. So, $\rho_{O}$ is nonnegative. Obviously, $\rho_{O}$ is symmetric.

Let $\mu_1= \mu_2 = \mu$. There exists a monetary risk measure $\mu_{1\, 2}\in O(\Delta(X))$ such that $\mu_{1\, 2}\in O(\pi^{1\,2}_{1})^{-1}(\mu_{1})\cap O(\pi^{1\,2}_{2})^{-1}(\mu_{2})$. Then
\begin{gather*}
\rho_{O}(\mu_1,\, \mu_2) \le \max\{\rho(x,\, x):\, (x,\, x)\in \mbox{supp}\,\mu_{1\, 2}\} = 0.
\end{gather*}
Inversely, let $\rho_{O} (\mu_1,\, \mu_2) = 0$. Then there exists a $\mu_{1\, 2}\in O(\pi^{1\,2}_{1})^{-1}(\mu_{1})\cap O(\pi^{1\,2}_{2})^{-1}(\mu_{2})$ such that $\rho(x,\, y) = 0$ for all $(x,\, y)\in \mbox{supp}\, \mu_{1\, 2}$. Consequently, $\mbox{supp}\, \mu_{1\, 2}$ must lie in the diagonal $\Delta(X)$. We have $(\varphi\circ \pi^{1\, 2}_{1})|_{\Delta(X)} = (\varphi\circ \pi^{1\, 2}_{2})|_{\Delta(X)}$, $\varphi\in C(X)$. From here with respect to Lemma 4 \cite{Rad1998} (see Page \pageref{coinsidemeasures} of the present paper)
\begin{gather*}
\mu_{1}(\varphi) = \mu_{1\, 2}(\varphi\circ\pi^{1\, 2}_{1}) = \mu_{1\, 2}(\varphi\circ\pi^{1\, 2}_{2}) = \mu_{2}(\varphi),\qquad \varphi\in C(X),
\end{gather*}
i.~e. $\mu_{1} = \mu_{2}$.

Let us show that the triangle inequality is true as well. Take an arbitrary triple $\mu_i\in O(X)$, $i=1,\,2,\,3$. Let $\mu_{1\,2},\,\mu_{2\,3}\in O(X^2)$ be $(\mu_1,\, \mu_2)$- and $(\mu_2,\, \mu_3)$-admissible monetary risk measures such that $\rho_{O}(\mu_1,\, \mu_2) = \max\{\rho(x,\, y):\, (x,\, y)\in\mbox{supp}\, \mu_{1\,2}\}$ and $\rho_{O}(\mu_2,\, \mu_3) = \max\{\rho(x,\, y):\, (x,\, y)\in\mbox{supp}\, \mu_{2\,3}\}$, respectively. Using Corollary \ref{overmu2} one can point out a $\mu_{1\,2\,3}\in O(X_{1\,2\,3})$ which satisfies the equalities
\begin{gather*}
O(\pi^{1\,2\,3}_{1\,2})(\mu_{1\,2\,3}) = \mu_{1\,2}\qquad \mbox{ and }\qquad O(\pi^{1\,2\,3}_{2\,3})(\mu_{1\,2\,3}) = \mu_{2\,3}.
\end{gather*}
We assume that $\mu_{1\,3} = O(\pi^{1\,2\,3}_{1\,3})(\mu_{1\,2\,3})$. Then  $\mu_{1\,3}$ is a $(\mu_1,\, \mu_3)$-admissible monetary risk measure. We have
\begin{multline*}
\rho_{O}(\mu_1,\, \mu_2) + \rho_{O}(\mu_2,\, \mu_3) = \\
= \max\limits_{(x_{1},\, x_{2})\in {\scriptsize\mbox{supp}}\,\mu_{1\,2}}\rho(x_{1},\, x_{2}) +
\max\limits_{(x_{2},\, x_{3})\in {\scriptsize\mbox{supp}}\,\mu_{2\,3}}\rho(x_{2},\, x_{3}) =\\
= \max\limits_{(x_{1},\, x_{2},\, x_{3})\in {\scriptsize\mbox{supp}}\,\mu_{1\,2\,3}}\rho(x_{1},\, x_{2}) +
\max\limits_{(x_{1},\, x_{2},\, x_{3})\in {\scriptsize\mbox{supp}}\,\mu_{1\,2\,3}}\rho(x_{2},\, x_{3}) \ge\\
\ge \max\limits_{(x_{1},\, x_{2},\, x_{3})\in {\scriptsize\mbox{supp}}\,\mu_{1\,2\,3}}\{\rho(x_{1},\, x_{2}) + \rho(x_{2},\, x_{3})\} \ge\\
\ge\max\limits_{(x_{1},\, x_{2},\, x_{3})\in {\scriptsize\mbox{supp}}\,\mu_{1\,2\,3}}\rho(x_{1},\, x_{3}) =\\
=\max\limits_{(x_{1},\, x_{3})\in {\scriptsize\mbox{supp}}\,\mu_{1\,3}}\rho(x_{1},\, x_{3}) \ge \rho_{O}(\mu_1,\,\mu_3),
\end{multline*}
i.~e. $\rho_{O}(\mu_1,\,\mu_3) \le \rho_{O}(\mu_1,\, \mu_2) + \rho_{O}(\mu_2,\, \mu_3)$.

For every pair of the Dirac measures $\delta_x$, $\delta_y$,  $x,\, y\in X$, the uniqueness of $(\delta_x,\, \delta_y)$-admissible monetary risk measure $\delta_{(x,\,y)}\in O(X^{2})$ implies that
\begin{gather*}
\rho_{O}(\delta_x,\, \delta_y) = \max\{\delta_{(x,\, y)}(\rho),\, (x,\, y)\in \mbox{supp}\,\delta_{(x,\,y)}\} = \rho(x,\,y).
\end{gather*}
From here we get that $\rho_{O}$ is an extension of $\rho$.

\end{Proof}

The following affirmation states one of remarkable properties of the metric $\rho_{O}$.

\begin{Prop}
{\rm $\mbox{diam}(O(X),\, \rho_{O}) = \mbox{diam}(X,\, \rho)$.}
\end{Prop}
\begin{Proof}
The proof is obvious.

\end{Proof}

\begin{Th}
{\rm The metric $\rho_{O}$ generates point-wise convergence topology on $O(X)$.}
\end{Th}
\begin{Proof}
Let $\{\mu_n\} \subset O(X)$ be a sequence and $\mu_0 \in O(X)$. Suppose the sequence converges to $\mu_0$ with respect to the point-wise convergence topology but not by $\rho_{O}$. Passing in the case of need to a subsequence, it is possible to regard that
\begin{gather*}
\rho_{O}(\mu_n,\, \mu_0)\ge a> 0\qquad \mbox{ for all positive integer }\quad n.
\end{gather*}
Consider an open neighbourhood of the diagonal $\Delta(X)\subset X^{2}$:
\begin{gather*}
U = \left\{(x,\, y)\in X^2:\, \rho(x,\, y)<\frac{a}{2}\right\}.
\end{gather*}
By virtue of Proposition \ref{diagneigh} there exist a positive integer $n$ and a $(\mu_0,\, \mu_n)$-admissible measure $\mu_{0\, n}\in O(X^2)$ such that $\mbox{supp}\,\mu_{0\, n}\subset U$, and consequently
\begin{gather*}
\rho_{O}(\mu_n,\, \mu_0) \le\frac{a}{2} < a.
\end{gather*}
The obtained contradiction finishes the proof.
\end{Proof}

\end{document}